\documentclass[12pt]{amsart}

\newtheorem{theorem}{Theorem}[section]
\newtheorem{proposition}{Proposition}[section]

\def\C{{\mbox{\rm\kern.24em
\vrule width.03em height1.43ex depth-.052ex \kern-.26em C}}}
\def\QSet{\mbox{\rm\kern.24em
\vrule width.03em height1.48ex depth-.051ex \kern-.26em Q}}
\def\Z{{\mbox{\rm\kern.25em
\vrule width.03em height0.57ex depth0ex
\kern.033em
\vrule width.03em height1.52ex depth-0.96ex \kern-.338em Z}}}

\def\N{{{\mbox{\rm I\kern-.2em N}}_0}}
\def\M{{{\mbox{\rm I\kern-.2em M}}}}
\def\PSet{\mbox{\rm I\kern-.22em P}}
\def\R{{\mbox{\rm I\kern-.22em R}}}

\def\I{{\bf I}}

\def\T{{\bf T}}

\def\I{{\bf I}}

\def\<{\left<}
\def\>{\right>}

\title[Variational estimates for paraproducts]
{Variational estimates for paraproducts}

\author
[Y. Do, C. Muscalu, C. Thiele]
{Yen Do \ \ \ Camil Muscalu \ \ \ Christoph Thiele }

\address{Yen Do, School of Mathematics,
Georgia Institute of Technology, Atlanta, GA 30332-0160, USA}
\email{yendo@math.gatech.edu}

\address{Camil Muscalu, Department of Mathematics,
Cornell University, Ithaca, NY 14852-4201 , USA}
\email{camil@math.cornell.edu}

\address{Christoph Thiele, Department of Mathematics,
UCLA, Los Angeles, CA 90095-1555, USA}
\email{thiele@math.ucla.edu}


\thanks{This material is based upon work supported by the National Science Foundation under grant DMS-0635607002. Any opinions, findings, and conclusions or recommendations expressed in this material are those of the author(s) and do not necessarily reflect the views of the National Science Foundation.}
\thanks{C.M. partially supported by NSF grant DMS-0653519.} 
\thanks{C.Th. partially supported by NSF grant DMS 1001535.}
\date{\today}

\begin{document}

\begin{abstract}

We generalize a family of variation norm estimates
(Theorem~\ref{lineartheorem}) of L\'epingle with endpoint
estimates of Bourgain and Pisier-Xu to a family of variational estimates for 
paraproducts (Theorem~\ref{bileplem}), both in the discrete and the 
continuous setting. This expands on work of Friz and Victoir, our focus being
on the continuous case and an expanded range of variation exponents.
\end{abstract}

\maketitle

\section{Introduction and Main Theorem}

In this paper, a band limited function on $\R$ with band width $N$ has Fourier
transform supported in $\{2^{-N}<|\xi|<2^{N}\}$.
A collection of functions $(f_i)_{i\in \Z}$ such that the dilates
$f_i(2^{i}x)$
are band limited with uniform band width
is called a continuous Littlewood-Paley family.
The collection is called a discrete Littlewood-Paley family if 
$f_i$ is spanned by the Haar functions associated to dyadic intervals of 
length $2^{1-i}$, it is then essentially a martingale. The following theorem is known:

\begin{theorem}\label{lineartheorem}

Given $1<p< \infty$, $1 \le r\le 2$, $r\le t<\infty$, and $N>0$, 
there is a constant $C$ such that for any collection $(N_k)_ {k\in \N}$ of measurable 
functions on $\R$ and any  Littlewood Paley data  $f_i$ (continuous with bandwidth $N$ 
or discrete) the following holds:

If $r< 2$ or $t > r$, then 
\begin{equation}\label{strong}
\|(\sum_{k} |\sum_{N_{k-1}<i\le N_k}
f_i |^t)^{1/t}\|_p\le C
\|(f_i)\|_{p,r} \ \ .
\end{equation}

If $r=2$ and $t=r$, then for every $\lambda>0$
\begin{equation}\label{weak}
\|\lambda ( \sharp \{k:|\sum_{N_{k-1}<i\le N_k}
f_i |>\lambda\})^{1/t} \|_p
\le C \|(f_i)\|_{p,r} \ \ .
\end{equation}
Here we have set
\begin{equation}\label{mixednorm}
\|(f_i)\|_{p,r}
:=\|(\sum_i |f_i|^r)^{1/r}\|_{p}\ \ .
\end{equation}
\end{theorem}

In the discrete case, these inequalities are special cases of known
inequalities for martingales: (\ref{strong}) is L\'epingle's inequality \cite{lepingle}
if $t>r$  and Pisier and Xu's inequality \cite{pisierxu} if $t=r$.
Inequality (\ref{weak}) was used by Bourgain \cite{bourgain} to give a proof of L\'epingle's 
inequality. Inequality (\ref{weak}) for $r<2$ or $t>r$ follows from (\ref{strong})
by Chebysheff's inequality.

For $r=2$, the continuous case of Theorem~\ref{lineartheorem} has been
observed by several authors (\cite{bourgain, jonesetal, 
jonesseegerwright}) in the following setting.
Define for some fixed function $f$
\begin{equation}\label{scalekernel}
f_i:=\phi_i*f, \ \ \phi_i(x):=2^{i}\phi(2^{i}x)\ \ ,
\end{equation}
where $\phi$ is some band limited test function. 
Then (\ref{mixednorm}) for $r=2$ is the norm of a Littlewood Paley square function and
dominated by $\|f\|_p$, 
and the domination is an equivalence if $\phi_i$ is chosen appropriately. 
We present a proof of the continuous case $r<2$ of Theorem~\ref{lineartheorem}
in Section \ref{contpx} as model for the unfortunately somewhat technical Section \ref{raiselowersect}.

We prove the following bilinear variant of Theorem:
\ref{lineartheorem}
\begin{theorem}\label{bileplem}
Given $1<p,q<\infty$, $1 \le r,s \le 2$,
$1/t\le 1/r+1/s$, $N>0$
there is a constant  $C$ such that for any collection $(N_k)$ of measurable functions
and any Littlewood Paley data $f_i$ and $g_j$ (each continuous with bandwidth $N$ or discrete) 
the following holds:

If $\max(r,s)<2$ or $t>rs/(r+s)$, then
\begin{equation}\label{bistrong}
\|(\sum_{k} |\sum_{N_{k-1}< i<j \le N_k} f_i\,   g_j |^t)^{1/t} \|_{\frac{pq}{p+q}} 
\le C 
\|(f_i)\|_{p,r}
\|(g_j)\|_{q,s} \ \ .
\end{equation}
If $\max(r,s)=2$ and $t=rs/(r+s)$, then
\begin{equation}\label{biweak}
\|\lambda (\sharp \{k: |\sum_{N_{k-1}< i<j \le N_k} f_i\,   g_j |>\lambda\})^{1/t} \|_{\frac{pq}{p+q}} 
\le C 
\|(f_i)\|_{p,r}\|(g_j)\|_{q,s} \ \ .
\end{equation}

\end{theorem}
For $r=s=2$ and $p=q$, $f_i=g_i$, the discrete case of (\ref{bistrong})
is a special case of a martingale inequality of Friz and Victoir \cite{frizvictoir}. 
All other instances of this theorem appear to be new.
Inequality (\ref{biweak}) for $\max(r,s)<2$ or $t>\frac{rs}{r+s}$ follows from (\ref{bistrong})
by Chebysheff. Note that Theorem \ref{bileplem} has a continuous, a discrete, and a mixed continuous and
discrete case, though the latter is maybe less natural.

Theorems \ref{lineartheorem} and \ref{bileplem} have alternative formulations, 
where $f_i$ and $g_j$ are not required to be band limited. Instead, in the continuous 
case we replace $f_i$ and $g_j$ on the left-hand-side by $\phi_i*f_i$ and
$\phi_j*g_j$ with $\phi_i$  as in (\ref{scalekernel}), while in the discrete
case we use Haar projections on the left-hand-side. In particular, we have
\begin{theorem}\label{nonband}
Let $\phi_i$ be as in (\ref{scalekernel}).
Given $1<p,q<\infty$, $1 \le r,s \le 2$,
$1/t\le 1/r+1/s$,
there is a constant  $C$ such that for any collection $(N_k)$ of measurable
functions and any sequences of functions
$f=(f_i)$, $g=(g_j)$ with the normalization
 \begin{equation}
\label{pqnorm}
\|(f_i)\|_{p,r}=\|(g_j)\|_{q,s} =1
\end{equation}
the following holds:

If $\max(r,s)<2$ or $t>rs/(r+s)$, then
\begin{equation}\label{bistrongnon}
\|(\sum_{k} |\sum_{N_{k-1}< i<j \le N_k} (\phi_i*f_i)\, (\phi_j*  g_j) |^t)^{1/t} \|_{\frac{pq}{p+q}} 
\le C  \ \ .
\end{equation}
If $\max(r,s)=2$ and $t=rs/(r+s)$, then
\begin{equation}\label{biweaknon}
\|\lambda (\sharp \{k: |\sum_{N_{k-1}< i<j \le N_k} (\phi_i*f_i)\,   (\phi_j*g_j) |>\lambda\})^{1/t} \|_{\frac{pq}{p+q}} 
\le C  \ \ .
\end{equation}
\end{theorem}
Theorem (\ref{bileplem}) implies Theorem \ref{nonband} by an application of the former
to the Littlewood Paley family
$(\phi_i*f_i)$ and by the Fefferman-Stein inequality \cite{stein} for the
Hardy Littlewood maximal operator $M$: 
$$\|(\phi_i*f_i)\|_{p,r}\le C 
\|(M f_i)\|_{p,r}\le C
\|(f_i)\|_{p,r}\ \ .$$
Conversely, the continuous case of Theorem \ref{bileplem} follows by specializing Theorem (\ref{nonband})
to a continuous Littlewood Paley family $(f_i)$ and choosing $\phi$ 
with large band width such that $\phi_i*f_i=f_i$.
We will use both
formulations, Theorem \ref{bileplem} has simpler notation while
Theorem \ref{nonband} works better with truncations of $f$ and $g$
and interpolation.

Another alternative formulation arises from setting 
$$f_i:=\sum_{|I|=2^{1-i}} a_I \phi_I$$
where $a_I$ are coefficients and for each dyadic interval
$I$ the function $\phi_I$ is an $L^\infty$ normalized band limited bump
function adapted to $I$ in the sense of \cite{stein}. In particular it satisfies for $\epsilon>0$
$$|\phi_I(x)|\le C (1+(x-c(I))/|I|)^{-(1+\epsilon)}\ \ ,$$
which allows to estimate (\ref{mixednorm}) by 
$$C 
\|(M \sum_{|I|=2^{1-i}} a_I 1_I )\|_{p,r}\le C \|(\sum_{x\in I} |a_I|^r)^{1/r}\|_{L^p(x)}\ .$$
The $L^\infty$ normalization of $\phi_I$ makes $a_I$ have the same normalization as the values of 
the corresponding $f_i$.

The restriction $i<j$ in Theorem~\ref{bileplem} is analogous to the paraproduct
\cite{coifmanmeyer}
of two functions $f$ and $g$:
$$P(f,g)=\sum_{i<j} (\phi_i*f)(\phi_j *g)\ \ .$$
The endpoint of Theorem~\ref{bileplem} for $t=\infty$ is a classical maximal paraproduct inequality,
for example in the convolutional case we have
\begin{proposition}
For $1<p,q<\infty$,
\begin{equation}\label{maxparaproduct}
\|\sup_{N_0,N_1} |\sum_{N_0<i<j\le N_1} (\phi_i*f)(\phi_j *g)|\|_{pq/(p+q)}\le C_{p,q}\|f\|_p\|g\|_q\ \ .
\end{equation}
\end{proposition}

Inequality (\ref{strong}) has a trivial endpoint for $r=t=1$ 
by the triangle inequality. Similarly, the endpoint of (\ref{bistrong}) for $r=1$ or $s=1$ is easy and 
elaborated in detail in Proposition \ref{1prop}. This endpoint is proved by reduction to Theorem \ref{lineartheorem}.
and used as interpolation endpoint to prove part of Theorem \ref{bileplem}.
Note that while classical variation norms as in (\ref{strong}) are somewhat pointless for $t\le 1$, because 
then the choice of consecutive points $N_{k+1}=N_k+1$ is extremal, the bilinear variational
expression in (\ref{bistrong}) is meaningful for $t>1/2$.

We make two simple observations on Theorem~\ref{bileplem}. 
The variant of (\ref{bistrong}) without the paraproduct restriction $i<j$ follows
from Theorem~\ref{lineartheorem} by H\"older's inequality.
The high exponent case $t>2$ of Theorem~\ref{bileplem}
has a simple proof that we sketch in the discrete case.  Let $\Delta_m$ denote the projection onto 
the Haar functions associated to dyadic intervals of length $2^{1-m}$. Then 
$$\Delta_m(\sum_{i<j} f_i g_j)
= \sum_{i<m} f_i g_m \ \ ,$$
$$\sum_{N_{k-1}< i<j \le N_k} f_i g_j
=\sum_{N_{k-1}<m\le N_{k}}\Delta_m[ \sum_{i<j}f_i g_j]
- \sum_{i\le N_{k-1}}f_i \sum_{N_{k-1}<j\le N_{k}} g_j\ \ .$$
The first term on the right hand side can be estimated by (\ref{strong})
and (\ref{maxparaproduct}), while the second term can be estimated by the Hardy Littlewood
maximal theorem and (\ref{strong}) applied to the terms 
terms involving $f$ and $g$ respectively.

We became interested in variational estimates for paraproducts
while studying $L^p$ estimates for a variational expression of the form
\begin{equation}\label{varbiest}
(\sum_{k}|\int_{N_{k-1}<\xi<\eta<N_k}
\widehat{f}(\xi)\widehat{g}(\eta) 
e^{i\pi x (\xi+ \eta)}\, d\xi d\eta|^r)^{1/r}\ \ .
\end{equation}
This can be viewed as a bilinear analogue of the variation norm Carleson operator
studied in \cite{osttw} or as variational variant of the bi-est operator introduced 
in \cite{MTT3}. It requires time-frequency analysis to be understood, and
Theorem~\ref{bileplem} provides the related lacunary
estimates. 

The general type of bilinear variational estimate as in 
Theorem~\ref{bileplem} and (\ref{varbiest}) is motivated by 
Terry Lyons' \cite{lyons} theory of ODE with rough driving terms.
This theory bootstraps Theorem~\ref{bileplem} 
to the diagonal case $p_1=\dots=p_m$ of the following multilinear generalization
of Theorem~\ref{bileplem} 
\begin{proposition}\label{multileplem}
For $2/M<r<\infty$, $1<p_m<\infty$,
and  $f_m\in L^{p_m}(\R^n)$ for $m=1,\dots,M$, we have 
with $1/p=\sum_m 1/p_m$:
$$
\|\sup_{N_k} (\sum_{k} |\sum_{N_{k-1}< i_1<\dots<i_M \le N_k} 
\prod_m (\phi_{i_m}* f_m)|^r)^{1/r} \|_{p} \le C \prod_m\|f_m\|_{p_m} 
\ .
$$
\end{proposition}
The non-diagonal case is likely to follow from a multilinear generalization of arguments as in 
this paper, though the exact elaboration of the details 
is beyond the scope of this paper.

Further consequences in the case of martingales are discussed in \cite{frizvictoir}.
The authors would like to thank Terry Lyons for pointing out Reference \cite{frizvictoir}
after circulation of an earlier draft of this paper.

\section{Proof of Pisier and Xu's inequality, continuous case}
\label{contpx}

Recall that Pisier and Xu's inequality is the case $t=r<2$ of inequality $(\ref{strong})$.
In the diagonal case, $p=t$, inequalities (\ref{strong}) and (\ref{weak}) can be written
as strong and weak type $L^p$ estimates (see \cite{stein} for this terminology)  
for mappings from functions on the measure space $\R\times \Z$ to functions
on the measure space  $\R\times \N$:
\begin{equation}\label{pxdiag}
\| \sum_{N_{k-1}(x)<i\le N_k(x)}
\phi_i*f_i(x) \|_{L^p(x,k)}
\le C
\|f_i(x)\|_{L^p(x,i)}\ \ ,
\end{equation}
\begin{equation}\label{doobdiag}
|\{(x,k):  |\sum_{N_{k-1}(x)<i\le N_k(x)}
\phi_i*f_i(x)|>\lambda\}|\le C\lambda^{-2}\|f_i(x)\|^2_{L^2(x,i)}\ \ .
\end{equation}
Here we use the analogue setup as in Theorem \ref{nonband}.
For a proof of (\ref{doobdiag}) we refer to \cite{jonesseegerwright}.
The endpoint of (\ref{pxdiag}) for $p=1$ follows from the
triangle inequality. Hence (\ref{pxdiag}) follows by Marcinkiewicz interpolation
from (\ref{doobdiag}), which completes the discussion of the diagonal case.

Starting from this diagonal case, we shall lower and raise the exponent $p$ by Calder\'on-Zygmund decomposition and sharp function techniques respectively. For the sub-diagonal case, $p<t$, fix $r$ and
employ Marcinkiewicz interpolation for $l^r$- valued functions.
The interpolation endpoints are the diagonal case $p=r$
and the weak type bound 
\begin{equation}\label{wt}
|\{x: (\sum_{k} |\sum_{N_{k-1}<i\le N_k}
\phi_i* f_i(x) |^r)^{1/r}>\lambda \}| \le \frac{C}{\lambda}
\|(\sum_i |f_i|^r)^{1/r}\|_{1}
\end{equation}
at $L^1$.
The latter follows via a Calder\'on-Zygmund decomposition of the vector function
$f=(f_i)$ at level $\lambda$. 
Let
$$E=\{x: M((\sum_i |f_i|^r)^{1/r})> \lambda\}$$
and let $\I$ be the collection of maximal dyadic intervals
contained in the set $E$.
Split $f$ as 
$$f=g+b=g+\sum_I b_{I}\ \ ,$$
where on each interval $I$ the function $g$ is constant equal to the average of $f$ on $I$, 
and each $b_I$ is $f-g$ restricted to $I$.
It suffices to prove (\ref{wt}) separately with $f$ replaced by
$g$ and $b$ on the left-hand-side. But
$$\|\phi_i*b_{i,I}\|_{L^1((3I)^c)}\le C\lambda |I|\min(2^i |I|, (2^{i}|I|)^{-\epsilon})$$
by smoothness $(|I|<2^{-i}$) and decay $(|I|>2^{-i})$ estimates for $\phi_i$. Hence, by
embedding $l^1$ into $l^r$,
$$\|
(\sum_{k} |\sum_{N_{k-1}<i\le N_k}
\phi_i* \sum_I b_{i,I} |^r)^{1/r}\|_{L^1((\bigcup 3I)^c)} $$
$$\le C \sum_I \sum_i \| 
\phi_i* b_{i,I}  \|_{L^1(( 3I)^c)}\le C\lambda |E| \le C\|f\|_{1,r}\ \ .$$
Then (\ref{wt}) for $b$ follows from Chebysheff's inequality. 
On the other hand,
$$\|g\|^r_{r,r}\le C\lambda^{r-1}\|g\|_{1,r}\le C\lambda^{r-1}\|f\|_{1,r}$$
and (\ref{wt}) for $g$ follows from the known diagonal estimate
and Chebysheff's inequality. This completes the proof of the
sub-diagonal case.

For the super-diagonal case, $p>r$, consider the 
sharp function 
\begin{equation}\label{sharpfunction}
(Tf)^{\sharp }(x)=\sup_{x\in I}\inf_c \frac{1}{|I|}\int_I |Tf(y)-c|\, dy\ \ ,
\end{equation} 
$$Tf(x)=
\sup_{(N_k)}(\sum_k|
\sum_{N_{k-1}<i\le N_k}
\phi_i*f_i(x) |^r)^{1/r} \ \ .$$
Here the supremum is over all sequences $(N_k)$. 
The $L^p$ norms of $Tf$ and $(Tf)^{\sharp}$ are comparable \cite{stein}, 
hence it suffices to estimate pointwise
\begin{equation}\label{sharpestimate}
(Tf)^{\sharp }(x) \le C   M_r ((\sum_i |f_i|^r)^{1/r}))(x)\ \ ,
\end{equation}
where
$$M_rh(x):=(\sup_{x\in I}\frac 1{|I|}\int_I |h(y)|^r\, dy)^{1/r}.$$

Fix an interval $I$ and let
$\tilde{f}_i$ be a constant function on $\R$ whose
value equals the average of $\phi_i*f_i$ on $I$
if $2^{-i} > |I|$ and equals zero
if $2^{-i}\le |I|$. Instead of taking infimum in $c$, we evaluate the definition
of the sharp function $(Tf)^{\sharp}$ with
$$c=\sup_{(N_k)}(\sum_k|
\sum_{N_{k-1}<i\le N_k}
\tilde{f}_i |^r)^{1/r}\ \ . $$
By the general norm inequality $|\|a\|-\|b\||\le \|a-b\|$
we have
$$|\sup_{(N_k)}(\sum_k|
\sum_{N_{k-1}<i\le N_k}
\phi_i*f_i(y) |^r)^{1/r}-c|$$
$$\le 
\sup_{(N_k)}(\sum_k|
\sum_{N_{k-1}<i\le N_k}
\phi_i*f_i(y)-\tilde{f}_i(y) |^r)^{1/r}\ .$$
We write
$$\phi_i*f_i-\tilde{f}_i= h_i^{[1]}+h_i^{[2]}+h_i^{[3]}\ \ .$$
where 
$$h_i^{[2]}=h_i^{[3]}=0$$ 
for $2^{-i}>|I|$, while for $2^{-i}\le |I|$ we have
$h_i^{[1]}=0$ and 
$$h_i^{[2]}=\phi_i*(f_i 1_{3I})
\ \ , h_i^{[3]}=\phi_i*(f_i 1_{(3I)^c})\ \ .$$ 
Estimating the three summands separately, we
have for $h^{[1]}$, using embedding of $l^1$ into $l^r$ and
smoothness of $\phi_i*f_i$,
$$\frac1{|I|}\int_I 
\sup_{(N_k)}(\sum_k|
\sum_{N_{k-1}<i\le N_k}
h_i^{[1]}(y) |^r)^{1/r} \, dy$$
$$\le \frac1{|I|}\int_I \sum_{2^{-i} > |I|}|h_i^{[1]}(y) |\, dy
\le C \sum_{2^{-i} > |I|}(2^i|I|)Mf_i(x)\ \ .
$$
This is dominated by the right hand side of (\ref{sharpestimate}).
By the diagonal estimate, we have
$$\frac{1}{|I|}\int_I \sup_{(N_k)}(\sum_k |
\sum_{N_{k-1} < i \le N_k}
|h_i^{[2]}(y)|^r)^{1/r})\, dy$$
$$\le C |I|^{-1/r}\|(\sum_{i}|f_i1_{3I}|^r)^{1/r}\|_r\ \ ,$$
which is again bounded by the right hand side of (\ref{sharpestimate}).
Finally, we have by decay estimates for $\phi_i$ and embedding
of $l^1$ into $l^r$,
$$\frac{1}{|I|}\int_I \sup_{(N_k)}
(\sum_{N_{k-1} < i \le N_k}
|h_i^{[3]}
(y)|^r)^{1/r}\, dy$$
$$\le \frac{1}{|I|}\int_I 
\sum_{2^{-i} \le |I|}
|h_i^{[3]}(y)|\, dy
\le 
C \sum_{2^{-i} \le |I|}(2^i|I|)^{-\epsilon} Mf_i(x)\ \ ,$$
which is again bounded by the right hand side of (\ref{sharpestimate}).
This concludes the proof of the super-diagonal case.

\section{Reduction of Theorem~\ref{bileplem} to the case $t=rs/(r+s)$}\label{bourgain-argument}

If $\max(r,s)<2$ we simply observe that (\ref{bistrong}) is monotone
in $t$ to reduce to the case $t=rs/(r+s)$. If $\max(r,s)=2$ we need to
show that the relatively weaker estimate (\ref{biweak}) for some $t$ implies the 
relatively stronger estimate (\ref{bistrong}) for $t_0>t$. This will be a bilinear 
variant of the argument used by Bourgain \cite{bourgain}.

We apply Marcinkiewicz interpolation in the following form. 
\begin{proposition}\label{biinter}
Let $f$, $g$ denote measurable functions on measure spaces $X,Y$ with values
in Banach spaces $F$ and $G$ respectively. Suppose we have
a bi-quasi-sublinear operation $(f,g)\to T(f,g)$ producing some measurable
function $T(f,g)$ on a measure space $Z$. Here we mean quasi-sublinear in each argument,
which is
$$T (f+\tilde{f},g)\le 
C (T(f,g) + T(\tilde{f},g))\ \ ,$$
$$T (f,g+\tilde{g})\le 
C (T(f,g) + T (f,\tilde{g}))\ \ .$$
Assume we have for all $\lambda>0$ 
the weak type estimates
$$\lambda |\{z:|T(f,g)(z)|>\lambda\}|^{1/p+1/q}\le C\|f\|_p \|g\|_q$$ 
for all corners $(p,q)$ of an axis parallel rectangle whose interior contains the point $(p_0,q_0)$, then 
we have the strong type estimate
$$\|T(f,g)\|_{p_0q_0/(p_0+q_0)} \le \tilde{C}\|f\|_{p_0} \|g\|_{q_0}\ \ .$$ 
\end{proposition}
The proof of this proposition follows the standard Marcinkiewicz argument, in the multilinear 
setting described in Janson \cite{janson}. One  splits both functions $f$ and $g$ according to small 
and large values (according to some level $\lambda$ that is later integrated on) and estimates 
the four resulting terms of $T(f,g)$ using the assumed estimates. While Janson requires sublinearity in each argument, 
the adaption to quasi-sublinearity is not difficult since we split $T(f,g)$ only into four terms
for each level $\lambda$.

Fix $p,q,r,s,t$ as in the theorem and let $t_0>t$. We shall work in the setting
of Theorem \ref{nonband} and deduce 
(\ref{bistrongnon}) for $t_0$ from (\ref{biweaknon}) for $t$.
By Proposition (\ref{biinter}) with $F=l^r(\Z)$ and $G=l^s(\Z)$ it
suffices to show that for any $\lambda>0$ and for any 
$f=(f_i),g=(g_j)$ with normalization (\ref{pqnorm})
we have
\begin{equation}\label{weakbourgain}
|\{x: (\sum_{k} F(k)^{t_0})^{1/t_0} 
> \lambda\}| \le C\lambda^{-\frac{pq}{p+q}} \ \ ,
\end{equation}
where
$$F(k) = |\sum_{N_{k-1}<i<j\le N_k} (\phi_i*f_i) (\phi_j*g_j)|\ \ .$$
Indeed, since the problem is invariant under dilation of the $x$ axis by powers of $2$,
it suffices to prove (\ref{weakbourgain}) for $\lambda=1$.
Let $$E = \{x: \sup_{k} F(k) > 1\}\ \ .$$ 
Then by (\ref{biweaknon}) for the exponent 
$t$ and with $u=pq/(p+q)$ we have
$$|E| \le \int |\{k: F(k) > 1\}|^{\frac ut} \le  C \ \ .$$
Considering level sets of $F(k)$, for $x\not\in E$ we have
$$\sum_{k} F(k)^{t_0}  
\le C \sum_{n\le 0} 2^{nt_0}  |\{k: F(k) >2^{n}\}|\ \ . $$
By H\"older together with a geometric sum if $u/t>1$ and by monotonicity of $l^v(\Z)$ in $v$ if $u/t\le 1$ we obtain 
$$(\sum_{k} F(k)^{t_0})  ^{u/t} \le C\sum_{n\le 0} 2^{n(1-\epsilon)ut_0/t} |\{k: F(k) >2^{n}\}|^{u/t}$$
for some small $\epsilon>0$. Using Chebysheff we have
$$|\{x\not\in E: \sum_{k} F(k)^{t_0}  >1\}| \le 
\int_{\mathbb R\setminus E} (\sum_{k} F(k)^{t_0})^{u/t} $$
$$\le C\sum_{n\le 0} 2^{n(1-\epsilon)ut_0/t} 
\int |\{k: F(k) >2^{n}\}|^{u/t}\le C\sum_{n\le 0} 2^{n(1-\epsilon)ut_0/t}2^{-nu}\ . $$
In the last inequality we have applied (\ref{biweaknon}).
The right-hand-side is summable for sufficiently small $\epsilon$, since $t_0>t$.
This proves (\ref{weakbourgain}) and completes the reduction to the case $t=rs/(r+s)$
for Theorem \ref{nonband}. Adaptions for the discrete and mixed cases are not difficult.
We shall assume $t=rs/(r+s)$
throughout the rest of this paper.

\section{The endpoint at $r=1$ or $s=1$.} \label{1section}

If $r=s=1$, hence $t=1/2$, we have the following trivial observation for $0<p,q\le \infty$,
which holds for arbitrary sequences $(f_i)$ and $(g_j)$: 
\begin{equation}\label{11case}
\| (\sum_k |\sum_{N_{k-1}
< i<j \le N_k} f_i\,   g_j |^{1/2})^{2}\|_{pq/(p+q)} 
\end{equation}
$$\le \| \Big(\sum_k (\sum_{N_{k-1}< i \le N_k} |f_i| )^{1/2} 
(\sum_{N_{k-1}< j \le N_k} |g_j| )^{1/2}\Big)^{2} \|_{pq/p+q}$$
$$\le \| (\sum_i |f_i|)(\sum_j |g_j|)\|_{pq/(p+q)}\le \|(f_i)\|_{p,1}\|(g_j)\|_{q,1}
\ \ .$$
Consider just one of $r$, $s$ equal to $1$, by symmetry of the argument below we may assume $r=1$.
We then have the following result:
\begin{proposition}\label{1prop}
Let $0<p\le \infty$ and $1<q<\infty$ and assume $t=s/(1+s)$.
Let $N_k$ by measurable functions, $(f_i)$ any sequence of measurable functions 
and let $(g_j)$ be a Littlewood Paley family.

If $1<s<2$, then
\begin{equation}\label{1strong}
\|(\sum_{k} |\sum_{N_{k-1}< i<j \le N_k} f_i\,  g_j
 |^t)^{1/t} \|_{\frac{pq}{p+q}} 
\le C 
\|(f_i)\|_{p,1}
\|(g_j)\|_{q,s} \ \ .
\end{equation}
If $s=2$, then 
\begin{equation}\label{1weak}
\|\lambda (\sharp \{k: |\sum_{N_{k-1}< i<j \le N_k} f_i\,  g_j  |>\lambda\})^{1/t} \|_{\frac{pq}{p+q}} 
\le C 
\|(f_i)\|_{p,1}\|(g_j)\|_{q,s} \ \ .
\end{equation}
\end{proposition}

Inequality (\ref{1strong}) follows quickly from Theorem \ref{lineartheorem}. Note that
$$|\sum_{N_{k-1}< i<j \le N_k} f_i\,   g_j |
\le (\sum_{N_{k-1}< i \le N_k} |f_i|)(  \sup_{N_{k-1}<i\le N_k} |\sum_{i<j\le N_k}   g_j |) \ \ . $$
Then estimate the left-hand-side of (\ref{1strong}) with H\"older's inequality by
$$\|(f_i)\|_{p,1} \| \big(\sum_k \sup_{N_{k-1}<i\le N_k} |\sum_{i<j\le N_k}   g_j |^s\Big)^{1/s}\|_q \ \ . $$
Introducing for each $x$ a new sequence which consists of  the sequence $N_k(x)$ interlaced with 
extremal choices of $i$ where applicable, we can estimate the second factor by
Theorem \ref{lineartheorem} and conclude (\ref{1strong}).

It remains to prove (\ref{1weak}). It will suffice to prove the analogue inequality with $g_j$
replaced by $\phi_j*g_j$ on the left-hand-side and $g_j$ not necessarily bandlimited
as in Theorem (\ref{nonband}).
We will use Proposition \ref{biinter}. Since the quantity
\begin{equation}\label{nonquasi}
\lambda (\sharp \{k: |\sum_{N_{k-1}< i<j \le N_k} (f_i)\, (\phi_j*g_j) |>\lambda\})^{1/t} 
\end{equation}
is not quasi subadditive, we shall replace it by an equivalent quasi subadditive quantity.
For $\lambda>0$ define $\rho_\lambda(x)=\min(|x|,\lambda)$ and
note that $\rho_\lambda(x-y)$ defines a metric distance between $x$ and $y$.
Hence the operation $T_\lambda$ with
\begin{equation}\label{tlambda}
T_\lambda (f,g)=
\lambda^{-1} \Big(\sum_{k}  \rho_\lambda(\sum_{N_{k-1}< i<j \le N_k} (f_i)  (\phi_j*  g_j) )^{2t} \Big)^{1/t}
\end{equation}
is quasi subadditive. It dominates (\ref{nonquasi}) and hence it suffices to show (here $r=1$)
$$
\| T_\lambda(f,g) \|_{\frac{pq}{p+q}} 
\le C 
\|f\|_{p,r}\|g\|_{q,s} \ \ .
$$
By Proposition~\ref{biinter} it suffices to prove
\begin{equation}\label{tlmu}
|\{x: T_\lambda (f,g) >\mu\}|
\le C\mu ^{-\frac{pq}{p+q}}\|f\|_{p,r}\|g\|_{q,s} 
\end{equation}
with $C$ independent of $\lambda$ and $\mu$. 
Switching back to expression (\ref{nonquasi}) we will prove 
\begin{equation}\label{weak-type-counting-function}
\Big|\{\lambda \big|\{k: |\sum_{N_{k-1}<i<j\le N_k} (f_i)(\phi_j*g_j)\}| > \lambda\}\big|^{1/t} > \mu\}\Big|
\end{equation}
$$\le C\mu ^{-\frac{pq}{p+q}}(\|f\|_{p,r}\|g\|_{q,s})^{\frac{pq}{p+q}}  \ \ .$$
We argue that this is sufficient to conclude (\ref{tlmu}). Consider for fixed $x$ the level sets of the function
$$F(k)= |\sum_{N_{k-1}< i<j \le N_k}  (f_i)  (\phi_j*  g_j)| \ \ ,$$
i.e. the sets of $k$ for which $2^{n}\lambda \le F(k)< 2^{(n+1)}\lambda$ for $n\le 0$.
Using geometric decay of $2^{\epsilon n}$ for some small $\epsilon>0$ we obtain
$$T_\lambda(f,g) \le C \sum_{n\le 0} 2^{(2-\epsilon)n}\lambda |\{k: F(k)>2^n\lambda\}|^{1/t} \ \ .$$
Then we can estimate
$$|\{x: T_\lambda(f,g) > \mu\}| $$
$$\le \sum_{n\le 0} \Big|\Big\{x: 2^{(2-\epsilon)n}\lambda |\{k: F(k)>2^n\lambda\}|^{1/t} > c 2^{n\epsilon} \mu \Big\} \Big|$$
$$\le C\sum_{n\le 0} (2^{(2\epsilon-1)n}\mu)^{-\frac{pq}{p+q}}(\|f\|_{p,r}\|g\|_{q,s})^{\frac{pq}{p+q}}\ .$$
Taking $\epsilon<1/2$ we obtain (\ref{tlmu}).

Now we prove  (\ref{weak-type-counting-function}). 
By dilating the $x$-axis, we may assume that the right-hand-side of (\ref{weak-type-counting-function})
is $C$. By multiplying $f$, $g$, $\lambda$, and $\mu$ by
appropriate constants we may assume the normalization (\ref{pqnorm}) and by the previous also $\mu=1$.
We write
$$|\{k: |\sum_{N_{k-1}<i<j\le N_k} (f_i)(\phi_j*g_j)|>\lambda\}| \le$$
$$ |\{k: \sum_{N_{k-1}<i<N_k} |f_i|>\lambda^{t}\}|$$
$$ +  |\{k: \sup_{N_{k-1}<i< N_k}|\sum_{i<j\le N_k} \phi_j*g_j|>\lambda^{t/s}\}|\ \ .$$
Therefore the left-hand side of (\ref{weak-type-counting-function}) has the upper bound
$$ \Big|\Big\{x: \lambda\big|\{k: \sum_{N_{k-1}<i<N_k} |f_i|>\lambda^{t}\}\big|^{1/t}>c \Big\}\Big|$$
$$ +\Big|\Big\{x: \lambda\big|\{k: \sup_{N_{k-1}<i< N_k}|\sum_{i<j\le N_k} \phi_j*g_j|>\lambda^{t/s}\}\big|^{1/t}>c \Big\}\Big|\ \ .$$
Using Chebysheff twice, the first term can be estimated by
$$|\{x: \sum_i |f_i|>c^t\}|\le C \|(f_i)\|_{p,1}^p\ \ .$$
The second term can be estimated by
$$C\left\|\lambda^{t/s} \big|\{k: \sup_{N_{k-1}<i< N_k}|\sum_{i<j\le N_k} \phi_j*g_j|>\lambda^{t/s} \}\big|^{1/s}\right\|^q_q 
\le  C \|g\|^q_{q,s}$$
by Theorem~\ref{lineartheorem}, applied for a sequence $N_k$ interlaced with elements $i$.
This proves (\ref{weak-type-counting-function}) 
and completes the proof of Proposition \ref{1prop}.

\section{The diagonal case of Theorem~\ref{bileplem} and $L^2$ theory}\label{diagonal-section}

The goal of this section is to prove Theorems~\ref{bileplem} and \ref{nonband} for 
one pair of exponents $p,q$ for any given triple $r,s,t$
with $t=rs/(r+s)$. In Section \ref{raiselowersect}, we will raise and lower $p,q$ 
using a bilinear variant of the argument in Section~\ref{contpx}, and thereby 
complete the proof of Theorem \ref{bileplem}. As in the previous section, the 
desired estimates are technically easier under some diagonal assumptions on $p,q$. 
We will assume throughout this section $p=r$ and $q=s$ and in particular 
$
t=pq/(p+q)$,
which allow us to turn (\ref{biweaknon}) into a weak type formulation analogous to 
(\ref{doobdiag}), namely 
the estimate (formulation as in Theorem \ref{nonband})
\begin{equation}\label{weakdiag}
\int \sharp \{k: |\sum_{N_{k-1}< i<j \le N_k} (\phi_i*f_i)\,   (\phi_j*g_j) |>\lambda\}\, dx
\le C \lambda^{-t} 
\end{equation} 
with the normalization assumption (\ref{pqnorm}).

Note that for $r=1$ or $s=1$ we have already proven (\ref{weakdiag})
in the previous section for many values of $p,q$. By interpolation
it suffices to prove (\ref{weakdiag})
for $r=s=2$. Proposition \ref{biinter} applied to the square with
corners $(2,2)$, $(2,1)$, $(1,2)$, $(1,1)$ will prove
the diagonal case of (\ref{bistrongnon}), hence (\ref{weakdiag}), in the open interior of the square. Note that we apply
Proposition \ref{biinter} with $F,G=\C$ and $X,Y=\R \times \mathbb Z$.

On the edges of the square we only need to prove the weak type estimate (\ref{weakdiag}),
which does not require Marcinkiewicz interpolation but simple interpolation
by one time truncation of $f$ or $g$ and applying the endpoint estimates.

The rest of this section will be concerned with proving the case $r=s=p=q=2$
and $t=1$ of (\ref{weakdiag}).

\subsection{Discrete case, stopping times}

As first step, we assume here that $(N_k)$ is a stopping time, which means that 
if $I$ is the dyadic interval of length $2^{-N_k(x)}$ containing $x$, 
then for all $y\in I$ we have $N_k(x)=N_k(y)$.
Such an interval $I$ is called a tree top and we can form the 
tree $T$ (or stopping time region) consisting of all dyadic intervals 
contained in $I$ but not contained in any smaller tree top. 

For $x$ in a tree top $I$ of length $|I|=2^{-N_{k-1}(x)}$ we can write
\begin{equation}\label{treelin}
\sum_{N_{k-1}<i<j\le N_{k}} f_i(x)g_j(x)
=
 \sum_{i<j} \Delta_i \Pi_T f(x)
\Delta_j \Pi_T g (x)
\end{equation}
where $f=\sum_i f_i$, $g=\sum_j g_j$ and $\Pi_T$ denotes the projection
onto the space spanned by the Haar functions associated to dyadic intervals
in the tree $T$. Here we have used that $\Delta_i f(x)$ equals the evaluation
of the projection of $f$ onto a single Haar function depending on $i$ and $x$ 
and for a stopping time $(N_k)$ this Haar function corresponds to an interval 
in the tree $T$ if and only if $N_{k-1}<i\le N_{k}$.

In the case of stopping times, we prove strong type (\ref{bistrong})
directly, which of course implies weak type (\ref{weakdiag}). Inserting (\ref{treelin}) on the left-hand-side of (\ref{bistrong}) 
and interchanging the order of sum and integration we obtain
$$
\sum_{T\in \T}  \|\sum_{i<j} (\Delta_i \Pi_T f)
(\Delta_j \Pi_T g)\|\ \ .
$$
By the paraproduct estimate (\ref{maxparaproduct}) this is bounded by a constant times
$$
 \sum_{T\in \T}  \|\Pi_T f\|_2
\|\Pi_T g\|_2\ \ .
$$
Applying Cauchy Schwarz in the sum in $T$ and using 
orthogonality of the projections $\Pi_T$ this implies (\ref{bistrong}).

\subsection{Discrete case, arbitrary $(N_k)$} 

We will prove weak type estimate (\ref{weakdiag}) for
arbitrary $(N_k)$ by comparing
the set of $\lambda$ jumps for $(N_k)$ to the set of $\lambda/4$ jumps
of an adapted stopping time $(\tilde{N}_k)$. Fix measurable functions $N_k$, discrete 
Littlewood Paley families $f_i$ and $g_j$, and $\lambda>0$.
For each $x$ choose $\tilde{N}_0(x)$ sufficiently close to $-\infty$ that the maximal paraproduct restricted to $i<j<\tilde{N}_0(x)$ is pointwise bounded by $\lambda/4$, and the tree tops for $\tilde{N}_0$ partition the real axis. Define recursively $\tilde{N}_k(x)$ to be the first time past 
$\tilde{N}_{k-1}(x)$ such that one of the following two conditions is satisfied
$$|\sum_{ \tilde{N}_{k-1}(x)<i<j \le \tilde{N}_k(x)} 
f_i(x) g_j(x) |\ge \lambda/4\ \ ,$$
$$
\sup_{\tilde{N}_{k-1}(x)<m<\tilde{N}_{k}(x)}
|\sum_{ \tilde{N}_{k-1}(x)<i \le m} 
f_i(x)|
|\sum_{ m<j \le \tilde{N}_k(x)} 
g_j(x)|\ge \lambda/4\ \ .$$
If no such $\tilde{N}_k(x)$ exists, we set 
$\tilde{N}_k(x)=\tilde{N}_{k-1 }(x)+1$.
The established estimate for stopping times gives good control
for the first kind of jumps. The integral over the number
of the second kind of jumps is controlled by 
$$   \int \sharp \{k:
|\sup_{N_{k-1}(x)<m<N_{k}(x)}
|\sum_{ N_{k-1}(x)<i \le m} 
(\Delta_i f(x))|
\ge c\lambda^{1/2}\}\, dx$$
$$+   \int \sharp \{k:
|\sup_{N_{k-1}(x)<m<N_{k}(x)}
|\sum_{ m<j \le N_k(x)} 
(\Delta_j g(x))|\ge c\lambda ^{1/2}\}\, dx\ \ ,$$
where each term can be controlled via (\ref{weak})
by $C\lambda$.

It remains to show that for every instance
$$|\sum_{N_{k-1}(x)<i<j\le N_{k}(x)} f_i(x)g_j(x)|>\lambda$$
there is a $\tilde{k}$ satisfying
\begin{equation}\label{squeeze}
N_{k-1}(x) < \tilde{N}_{\tilde{k}}(x)\le N_{k}(x)\ \ .
\end{equation}
Assume the contrary that so such $\tilde k$ exists, then by choice of $\tilde{N}_{0}$ one must have $\tilde{N}_0(x) \le N_{k-1}(x)$. Let $\tilde{k}$  be the largest index such that
$\tilde{N}_{\tilde{k}-1}(x)\le N_{k-1}(x)$. Then $\tilde{N}_{\tilde{k}}(x)> N_k(x)$ and we have
$$|\sum_{N_{k-1}(x)< i<j\le N_k(x)} f_i(x) g_j(x)|
\le |\sum_{\tilde{N}_{\tilde{k}-1}(x)< i<j\le N_k(x)} f_i(x) g_j(x)|$$
$$+ |\sum_{\tilde{N}_{\tilde{k}-1}(x)< i<j\le N_{k-1}(x)} f_i(x) g_j(x)|$$
$$+ |\sum_{\tilde{N}_{\tilde{k}-1}(x)< i\le N_{k-1}(x)}f_i(x)  
\sum_{N_{k-1}(x)< j\le N_k(x)} g_j(x)|$$
Each term on the right hand side is
less than $\lambda/4$ by choice of $\tilde{N}_{\tilde{k}}$, while the left-hand side is larger than $\lambda$.
This contradiction proves (\ref{squeeze}).

\subsection{Mixed and continuous case} 

We compare the discrete case with the mixed and continuous case of 
(\ref{weakdiag}) using an inequality by Jones, Seeger and Wright (\cite{jonesseegerwright}), namely 
that the square function
$$Sf(x)=(\sum_i |{\mathbb E}_i f -\varphi_i*f|^2)^{1/2}$$
is bounded in $L^2$,
where ${\mathbb E}_i=\sum_{j\le i}\Delta_j$ and
$\widehat{\varphi_i}$ is supported in 
$|\xi|\le 2^N$ and constant $1$ on  $|\xi|\le 2^{N-1}$.

We write the telescopic sum 
$$
\sum_{N_{k-1}< i<j \le N_k} (\Delta_i f) (\Delta_j g)
=\sum_{N_{k-1}< j \le N_k} 
({\mathbb E}_{j-1}
 f-{\mathbb E}_{N_{k-1}}f)  (\Delta_j g)$$
and compare with 
$$
\sum_{N_{k-1}< j \le N_k} 
(f*\varphi_{j-1}- f*\varphi_{N_{k-1}})  (\Delta_j g)\ \ .$$
The difference can be estimated with another telescopic sum by
\begin{equation}\label{martcont}
|\sum_{N_{k-1}< j \le N_k} 
({\mathbb E}_{j-1}f - f*\varphi_{j-1})  (\Delta_j g)|
\end{equation}
$$+ | 
({\mathbb E}_{N_{k-1}}f - f*\varphi_{N_{k-1}}) 
({\mathbb E}_{N_{k}}g-
{\mathbb E}_{N_{k-1}}g)
|\ .$$
Considering the first term, we estimate
$$\lambda \sharp \{k: |\sum_{N_{k-1}< i<j \le N_k} 
({\mathbb E}_{j-1}f - f*\varphi_{j-1})  (\Delta_j g)
 |>\lambda\}$$
 $$ \le \sum_j |({\mathbb E}_{j-1}f - f*\varphi_{j-1}
 )  (\Delta_j g)|$$
$$\le 
(\sum_j |{\mathbb E}_jf - f*\varphi_j|^2)^{1/2}  
(\sum_j |\Delta_j g|^2)^{1/2} \ \ .$$
The $L^{pq/(p+q)}$ norm of the right-hand-side can be estimated using
H\"older's inequality and the classical square function  estimate
together with the
square function estimate of Jones, Seeger and Wright.
Turning to the second  term in (\ref{martcont}) we have
$$\lambda \sharp \{k: |
({\mathbb E}_{N_{k-1}}f - f*\varphi_{N_{k-1}}) 
({\mathbb E}_{N_{k}}g-
{\mathbb E}_{N_{k-1}}g)
 |>\lambda\}$$
$$\le \lambda \sharp \{k: |
{\mathbb E}_{N_{k-1}}f - f*\varphi_{N_{k-1}}|>\lambda^{1/2}\}$$
$$+
\lambda \sharp \{k: |
{\mathbb E}_{N_{k}}g - 
{\mathbb E}_{N_{k-1}}g
|>\lambda^{1/2} \}\ \ .$$
Then we estimate the first term by the square function of
Jones Seeger Wright and the second term using Theorem~\ref{lineartheorem}.
This proves (\ref{weakdiag}) for the mixed
convolutional and discrete case.

Using a very similar argument one can obtain (\ref{weakdiag}) for the mixed
discrete (in $f$) and convolutional (in $g$) case, and by combination
one can obtain the pure convolutional case.

\section{Completion of the proof of Theorem~\ref{bileplem}}\label{raiselowersect}

For each pair $r,s$ and $t=rs/(r+s)$ the previous section 
proves Theorem \ref{bileplem} for the specific
pair $(p,q)=(r,s)$. Beginning with this pair, we will lower and raise the exponents 
$p$ and $q$ by Calder\'on-Zygmund decomposition and sharp function techniques respectively, 
as in Section \ref{contpx}. This will complete the proof of the theorem.

We will consider operators $S(f,g)$ and have two different cases: 
\begin{enumerate}

\item 
\label{case1}
The quasi-bi-sublinear case
$$S(f,g)=
\Big(\sum_{k}  |\sum_{N_{k-1}< i<j \le N_k} (\phi_i*f_i)  (\phi_j*  g_j) |^{t} \Big)^{1/t}$$
for $t<1$. This case is needed to prove (\ref{bistrongnon}) for $\max(r,s)<2$.

\item
\label{case2}
 The equivalent-to-quasi-bi-sublinear case
$$S(f,g)=
\lambda (\sharp \{k: |\sum_{N_{k-1}< i<j \le N_k} f_i\,   g_j |>\lambda\})^{1/t} $$
for arbitrary $\lambda$ and $t\le 1$. This case is needed to prove
(\ref{biweaknon}) when $\max(r,s)=2$.
\end{enumerate}

The desired estimate now becomes an $L^{pq/(p+q)}$ estimate for $S(f,g)$. We will discuss the more difficult Case \ref{case2} in more detail, the other case is similar but
easier.

In this section we focus on the continuous case. The mixed or discrete case follows essentially the same arguments,
with possible simplifications at some places. We shall no further comment on the mixed and discrete case.

\subsection*{Lowering exponents}
Fix $r,s$ and first consider Case \ref{case2}. We work in the setting of Theorem \ref{nonband}.
Assuming (\ref{biweaknon}) holds for specific exponents $(p_0,q_0)$, we prove (\ref{biweaknon})
for all $p,q$ with $1 < p < p_0$, $1 < q < q_0$. More precisely, passing from $S(f,g)$ to the
quasi bi-sublinear
operation 
$$T_\lambda (f,g)=
\lambda^{-1} \Big(\sum_{k}  \rho_\lambda(\sum_{N_{k-1}< i<j \le N_k} (\phi_i*f_i)  (\phi_j*  g_j) )^{2t} \Big)^{1/t}$$
similarly as in Section \ref{1section} it suffices to prove 
$$
\| T_\lambda(f,g) \|_{\frac{pq}{p+q}} 
\le C 
\|f\|_{p,r}\|g\|_{q,s}\ \ .
$$
By Proposition \ref{biinter} it suffices to prove
for all $1 \le  p \le  p_0$, $1 \le  q \le q_0$ and all $\mu>0$ the weak type bound
$$
|\{x:  T_\lambda(f,g)\ge \mu\}|\le C \mu^{-pq/(p+q)}\|f\|_{p,r}\|g\|_{q,s} \ \ .
$$
As in  Section \ref{1section} it suffices to prove 
$$
\Big|\{\lambda \big|\{k: |\sum_{N_{k-1}<i<j\le N_k} (\phi_i*f_i)(\phi_j*g_j)\}| > \lambda\}\big|^{1/t} > \mu\}\Big| $$
$$\le C \mu^{-pq/(p+q)}\|f\|_{p,r}\|g\|_{q,s} \ \ .
$$

By symmetry of the argument, it suffices to fix $q=q_0$ and vary $p$. By simple interpolation it 
suffices to consider the extremal cases. The endpoint $p=p_0$ is true by assumption, so we may assume $p=1$. 
As in Section \ref{1section} we may assume (\ref{pqnorm}) and $\mu=1$.

Split $f=(f_i)$ into a Calder\'on-Zygmund decomposition 
at level $1$ as in Section \ref{contpx}. The good function
is in $L^{p_0}$ and we can apply the known estimate and Chebysheff.
The bad function satisfies outside the expanded bad set the estimate
\begin{equation}\label{farfield}
\|\sum_i|\phi_i*b_i|\|_{L^1((\bigcup 3I)^c)}\le C\|f\|_{1,r} \ \ .
\end{equation}
Applying Proposition \ref{1prop} to the sequences $((\phi_i*b_1)1_{(\bigcup 3I)^c})$ and $(g_j)$ we get
the desired estimate for the bad function.
This completes the proof of the weak type estimate for $p=1$.

In Case \ref{case1} we proceed similarly, except we work directly with the quasi sublinear
expression
$$T_\lambda(f,g)= 
\Big(\sum_{k}  |\sum_{N_{k-1}< i<j \le N_k} (\phi_i*f_i)  (\phi_j*  g_j) )|^t \Big)^{1/t}\ \ .$$

\subsection*{Raising the exponents in Case \ref{case2}}

In conjunction with the previous subsection, it suffices to prove
(\ref{biweaknon}) for every $p>r$ and $q>s$ satisfying the ``diagonal'' condition 
\begin{equation}\label{diagonal-relation}
\frac{p}{r}= \frac{q}{s} \ \ ,
\end{equation}
which includes pairs with arbitrarily large $p$ and $q$. 
Consider the function 
$$T_\lambda (f,g)(x)=
\sup_{(N_k)}(\sum_k \rho_\lambda(
\sum_{N_{k-1}<i<j\le N_k}
(\phi_i*f_i)(x) \ (\phi_j*g_j)(x))^{2t})^{\frac 1{2t}} \ \ .$$
Similarly as before, (\ref{biweaknon}) follows
from bounds 
$$\|\lambda^{-1}T_\lambda(f,g)^2\|_{pq/(p+q)}
\le C $$
under the normalization assumption (\ref{pqnorm}).
These are $L^{2pq/(p+q)}$ bounds for $T_\lambda(f,g)$ which we can prove by the corresponding bound for the
sharp function $T_\lambda(f,g)^{\sharp}$. 


 
The desired estimate now follows from bounds
\begin{equation}\label{pointwise}
\lambda^{-1}(T_\lambda(f,g)^{\sharp})^2
\end{equation}
$$\le C  
M_{r}((\sum_i|f_i|^r)^{1/r}) \ \  M_{s}((\sum_j|g_j|^s)^{1/s}) $$
$$+C
\sup_{x\in I} \Big(\frac {\lambda^{1/2} }{|I|}\int_{I}
\sup_{N_k} 
\Big|\{k:  |\sum_{N_{k-1}< i,j \le N_k} (\phi_i*{f}_i)  (\phi_j*  g_j-{\tilde{g}_{j,I}}) | >\lambda\}
\Big|^{\frac{1}{2t}}\Big)^{2}$$
pointwise at every $x$.
Here we have defined the function $\tilde{g}_{j,I}$
with respect to the interval $I$ as in Section \ref{contpx}.

We will first show how to estimate the right-hand-side of (\ref{pointwise}).
The first term satisfies the desired $L^{pq/(p+q)}$ bound by H\"older and the maximal theorem. 
We estimate the second term pointwise by a constant times 
$$M_{1/2}(
\sup_{N_k} \lambda  \Big|\{k:  |\sum_{N_{k-1}< i \le N_k} \phi_i*{f}_i  
 |>\lambda^{t/r}\} \Big|^{1/t})$$
$$+
\sup_I \Big(\frac 1{|I|} \int_I 
\sup_{N_k} \lambda^{1/2}  \Big|\{k:  |\sum_{N_{k-1}< j \le N_k} \phi_j*  g_j-{\tilde{g}_{j,I}} |>\lambda^{t/s}\} \Big|^{1/(2t)}\Big)^2\ \ .$$
Since $pq/(p+q)>1/2$, we can estimate the $L^{pq/(p+q)}$ norm of the first term using the maximal theorem adapted to $M_{1/2}$, by
$$ C\|
\sup_{N_k} \lambda  |\{k:  |\sum_{N_{k-1}< i \le N_k} \phi_i*{f}_i|>\lambda^{t/r}\}|^{1/t}\|_{\frac{pq}{p+q}} \ \ .$$
Observe that $\frac{rpq}{t(p+q)} = p$ thanks to (\ref{diagonal-relation}). So by (\ref{weak}) with $\lambda^{t/r}$ in place of $\lambda$ we can estimate the last display by
$$C\|(\sum_{i} |f_i|^r)^{1/r}\|_{p}^{r/t} = C \ \ .$$

To estimate the $g_j$- term, we decompose as in Section \ref{contpx}
$$\phi_j*g_j - \tilde{g}_{j,I} = k_j^{[1]}+ k_j^{[2]} + k_j^{[3]} \ \ .$$


To estimate the contribution of $k_j^{[2]}$, we'll use a weak-type variant of the argument in Section~\ref{contpx}. First, by adapting the sequence $N_k$ one can bound this contribution by
$$\sup_I \frac 1{|I|^2} \Big\|\sup_{N_k} \lambda^{t/s}  |\{k:  |\sum_{N_{k-1}< j \le N_k} \phi_j*(g_j1_{3I})|>\lambda^{t/s}\}|^{1/s}\Big\|_{L^{s/(2t)}(I)}^{s/t} \ \ ,$$
and since $t\ge 1/2$ we can use H\"older to increase $L^{s/(2t)}(I)$ to $L^{s}(I)$,
$$\le \sup_I \frac 1{|I|^{1/t}} \Big\|\sup_{N_k} \lambda^{t/s}  |\{k:  |\sum_{N_{k-1}< j \le N_k} \phi_j*(g_j1_{3I})|>\lambda^{t/s}\}|^{1/s}\Big\|_{L^{s}(I)}^{s/t} \ \ .$$
The point is that now we can use the diagonal case of (\ref{weak}) (with $\lambda^{t/s}$ in place of $\lambda$). This gives the following estimate for the last display:
$$\le C\sup_I \frac 1{|I|^{1/t}} \Big\|(\sum_{j} |g_j1_{3I})|^s)^{1/s}\Big\|_{L^{s}}^{s/t} \ \ \le C \Big(M(\sum_j |g_j|^s)\Big)^{1/t} \ \ .$$

Since $\frac{pq}{t(p+q)} = \frac{q}{s} > 1$, we can estimate the $L^{pq/(p+q)}$ norm of the last display using the maximal theorem, by
$$ \le C\|\sum_j |g_j|^s\|_{q/s}^{1/t}  = C\|g\|_{q,s}^{s/t} = C \ \ .$$

The arguments for estimating the contributions of $k_j^{[3]}$ and  $k_j^{[1]}$ will be similar to each other, and below we'll only show details for the latter.

The contribution of $k_j^{[1]}$ is estimated using decay in $j$ as in Section~\ref{contpx}, by
$$\le \sup_I \Big(\frac 1{|I|} \int_I (\sup_{N_k} \sum_k |\sum_{N_{k-1}<j\le N_k} k_j^{[1]} |^s)^{1/(2t)}\Big)^2$$
$$\le C\sup_I  \Big(\frac 1{|I|} \int_I (\sum_j |k_j^{[1]}|)^{s/(2t)}\Big)^2 \le C\sup_I  \Big(\frac 1{|I|} \int_I (\sup_j Mg_j)^{s/(2t)}\Big)^2$$
$$\le C M_{1/2}((\sup_j Mg_j)^{s/t})(x)$$
which satisfies the desired $L^{pq/(p+q)}$ bound by the maximal theorem and (\ref{diagonal-relation}). This concludes the bound of the right-hand-side of (\ref{pointwise}).

It remains to  prove the bound (\ref{pointwise}) for the sharp function. 
Fix an interval $I$ and define constant functions $\tilde{f}$ and 
$\tilde{g}$ as in Section \ref{contpx}. We'll use the constant
$$c := \sup_{(N_k)}(\sum_k \sum_{N_{k-1}<i<j\le N_k}
\tilde{f}_i \ \tilde{g}_j)^{2t})^{\frac 1{2t}}$$
for the sharp function estimate. By metric properties we obtain
$$|T_\lambda(f,g)-c|$$
\begin{equation}\label{polar}
\le T_\lambda(f-\tilde{f},g-\tilde{g})+
T_\lambda(f-\tilde{f},\tilde{g})
+T_\lambda(\tilde{f},g-\tilde{g})\ \ .
\end{equation}
Focusing on the first term in (\ref{polar}), we estimate the average 
$$\lambda^{-1} ( \frac1{|I|}\int_I T_\lambda(f-\tilde{f},g-\tilde{g})(x) \, dx)^2$$
using the equivalent estimate with

$$  \frac1{|I|^2} \Big\|
\sup_{N_k} \lambda \Big|\{k:  |\sum_{N_{k-1}< i<j \le N_k} (\phi_i*f_i-\tilde{f}_i)  (\phi_j*  g_j-\tilde{g_j}) |>\lambda\}\Big|^{1/t}
\Big\|_{L^{1/2}(I)}\ \ .$$
We'll show that the contribution of $ T_\lambda(f-\tilde{f},g-\tilde{g})$ is bounded by the product of two maximal functions on the right-hand side of (\ref{pointwise}). We split as in Section \ref{contpx}
$$\phi_i*f_i-\tilde{f}_i= h_i^{[1]}+h_i^{[2]}+h_i^{[3]}\ \ ,$$
$$\phi_j*g_j-\tilde{g}_j = k_j^{[1]}+k_j^{[2]}+k_j^{[3]}\ \ .$$
The estimate for the term involving $h^{[2]}$ and $k^{[2]}$ follows from H\"older (to increase $L^{1/2}$ to $L^{t}$) and the known estimate for the diagonal case 
$(p,q)=(r,s)$. For all other estimates we can use the $1$-variation bound of the error terms as in Section~\ref{contpx} in conjunction with the endpoint estimates of Proposition~\ref{1prop}. The key idea is to incorporate the restriction to the interval $I$ into the error terms \emph{before} applying this Proposition, because we want to preserve locality for future estimates by maximal functions. For instance, to estimate the combination $(h^{[2]}, k^{[3]})$ we first equivalently replace $k^{[3]}_j$ by $k^{[3]}_j 1_I$, then increase $L^{1/2}$ to $L^{r/(r+1)}$ and decrease $t$ to $\frac{r}{r+1}$, then estimate using Proposition~\ref{1prop}, by
$$\le C \frac 1{|I|^{(r+1)/r}} \|(f_i1_{3I})\|_{r,r} \|(k_j^{[3]}1_I)\|_{1,1} \le C M_r ((\sum_{i}|f_i|^r)^{1/r}) \sup_j Mg_j \ \ .$$

Turning to the second term in (\ref{polar}), note that by support of $\tilde{g}_j$ we may restrict attention
to $j<-\log_2|I|$. Thanks to the constraint $i<j$, only the term $h_i^{[1]}$ in the splitting of the function 
$\phi_i*f_i-\tilde{f}_i$ then appears in the summation.
While $\tilde{g_j}$ itself does not have geometric decay in $j$, the geometric
decay of $h_i^{[1]}$ in $i$ is sufficient to obtain geometric summability
in both $i$ and $j$ thanks to the constraint $i<j$. Specifically, using $t\ge 1/2$ and H\"older,
$$\lambda^{-1} \frac1{|I|^2} (\int_I  T_\lambda(h^{[1]}, \tilde g))^2 \le \Big(\frac1{|I|}\int_I\Big(\sum_{i<j<-\log_2|I|} |h_i^{[1]}| |\tilde g_j|\Big)^{1/2} \Big)^2 $$
$$\le C\Big(\frac1{|I|}\int_I\Big(\sum_{i<j<-\log_2|I|} (2^i|I|) \inf_{y\in I} M f_i(y) \inf_{y\in I} Mg_j(y) \Big)^{1/2} \Big)^2 $$
$$\le C\sup_i M f_i(x) \, \sup_j Mg_j(x)  \ \ .$$
This proves the desired estimate for the second term. Combining 
the first and second terms in (\ref{polar}) we also obtain bounds for the variant of
the second term with $\tilde{g}$ replaced by $g$.

Turning to the last term in (\ref{polar}), we may also replace $\tilde{f}$ with $f$. To apply the same
argument as for the second term, we need to replace the constraint $i<j$ by $i\ge j$,
which can be done via triangle inequality and estimates on the unconstrained expression
$$  \frac1{|I|^2} \Big\|
\sup_{N_k} \lambda  \Big|\{k:  |\sum_{N_{k-1}< i,j \le N_k} (\phi_i*{f}_i)  (\phi_j*  g_j-\tilde{g_j}) |>\lambda\}\Big|^{1/t}
\Big\|_{L^{1/2}(I)}\ \ .$$
This term can be estimated by the second term on the right-hand-side of (\ref{pointwise}).
This completes the proof of the pointwise bound for the sharp function.

\subsection*{Raising exponents for Case \ref{case1}} The argument in this section is similar to and simpler than last section and we'll only comment on the necessary changes. Consider a pair $(p,q)$ such that $r<p<\infty$, $s<q<\infty$. To show the desired $L^{\frac{pq}{p+q}}$ bound for
$$S(f,g)=
\Big(\sup_{N_k}\sum_{k}  |\sum_{N_{k-1}< i<j \le N_k} (\phi_i*f_i)  (\phi_j*  g_j) |^{t} \Big)^{1/t}$$
we'll show the corresponding bound for $T(f,g)^{\sharp}$, where $T(f,g) = S(f,g)^{1/2}$. Below we'll show the pointwise bound:
\begin{equation}\label{pointwisecase1}
T(f,g)^{\sharp} \le C \Big(M_{r}((\sum_i |f_i|^r)^{1/r}) M_{s}((\sum_j |g_j|^s)^{1/s})\Big)^{1/2}
\end{equation}
$$+ C \Big(  M(\|\phi_i*f_i\|_{V^r}) M_{s}((\sum_j|g_j|^s)^{1/s})\Big)^{1/2} \ \ .$$
Using this inequality, the desired $L^{\frac{2pq}{p+q}}$ bound for $T(f,g)^{\sharp}$ follows from H\"older's inequality, the maximal theorem and Theorem~\ref{lineartheorem}.

Below we show (\ref{pointwisecase1}). Fix any interval $I$ and define $\tilde{f}$ and $\tilde{g}$ as in Section~\ref{contpx}. We'll use the following constant for the sharp function estimate
$$c = \Big(\sup_{N_k}\sum_{k}  |\sum_{N_{k-1}< i<j \le N_k} \tilde{f}_i \ \tilde{g}_j |^{t} \Big)^{1/t} \ \ .$$
Since $\frac 12 \le t < 1$, both $|x|^t$ and $|x|^{1/(2t)}$ define a norm, so $T(f,g)$ is sublinear. We actually do not need $t<1$ for this to be true, because when $t\ge 1$ we could combine the $|x|^{1/2}$ and $\ell^t(\mathbb Z)$ norms. We then have
$$|T(f,g)-c| \le T(f-\tilde f, g-\tilde g) + T(f-\tilde f, \tilde g) + T(\tilde f, g-\tilde g) \ \ .$$ 

The contributions of $T(f-\tilde f, g-\tilde g)$ and of $T(f-\tilde f, \tilde g)$ can be estimated as in the consideration of Case~\ref{case2} in the last section. To estimate the contribution of $T(\tilde f, g- \tilde g)$, as before it suffices to estimate the unconstrained expression:
$$\frac 1{|I|} \int_I \sup_{N_k} (\sum_k |\sum_{N_{k-1}<i,j\le N_k} \tilde f_i (\phi_j*g_j - \tilde g_j)|^t)^{\frac 1{2t}} \ \ .$$
By H\"older's inequality we can estimate this by
$$\le \sup_{N_k}(\sum_k |\sum_{N_{k-1}<i\le N_k} \tilde f_i|^r)^{\frac 1{2r}} \frac 1{|I|} \int_I \Big(\sup_{N_k}\sum_k |\sum_{N_{k-1}<i\le N_k}(\phi_j*g_j - \tilde g_j)|^s\Big)^{\frac 1{2s}} \ \ .$$
Using the estimates in Section~\ref{contpx} we can estimate the last display by
$$\le C\Big(\sup_{N_k}(\sum_k |\sum_{N_{k-1}<i\le N_k} \tilde f_i|^r)^{\frac 1{r}}\Big)^{1/2} \Big(M_{s}((\sum_j|g_j|^s)^{1/s})\Big)^{1/2} \ \ .$$
Now, since $\tilde{f}_i = \frac 1{|I|}\int_I \phi_i * f_i$ for $i < -\log_2 |I|$ and $\tilde{f}_i=0$ otherwise, by Minkowski's inequality we can estimate the last display by
$$\le C\Big(\frac 1{|I|} \int_I(\sup_{N_k}\sum_k |\sum_{N_{k-1}<i\le N_k} \phi_i* f_i|^r)^{\frac 1{r}}\Big)^{1/2} \Big( M_{s}((\sum_j|g_j|^s)^{1/s})\Big)^{1/2}$$
$$\le C \Big(  M(\|\phi_i*f_i\|_{V^r}) M_{s}((\sum_j|g_j|^s)^{1/s})\Big)^{1/2} \ \ .$$

\end{document}